\documentclass[10pt]{amsart}
\usepackage{mystyle}

\title[Anisotropy via Geometry of Numbers]{On bond anisotropy for the Ising model through the geometry of numbers}
\author{Ren\'e R\"uhr}

\begin{document}

\begin{abstract}
  We associate to each unit volume lattice of $\R^d$ the Ising model with bond variables equal to the inverse successive minima of that lattice.
This induces the notion of a critical temperature for a random lattice for which integrability exponents are proven.
\end{abstract}

\maketitle
\setcounter{tocdepth}{1}

In this note we consider the classical Ising model \cite{ising1925beitrag} in dimensions $d\geq2$ associated to the Hamiltonian
\begin{equation}\label{eq:Hamiltonian_simple}
\mathcal{H}_{J, \beta, \Lambda}(\sigma) = -\beta\sum_{i\sim j} J_{ij}\sigma_i\sigma_j.
\end{equation}
The parameter $\beta=1/T$ denotes the inverse temperature. The sum ranges over all adjacent sites $i\sim j$.
Above a critical temperature $T_c(J)$, there is a unique Gibbs measure associated to $\mathcal{H}_{J, \beta, \Lambda}$.
The value $T_c(J)$ is known exactly for $d=2$ by a celebrated result of Onsager \cite{onsager1944crystal},
and for general $d$,
one has the mean field bounds of Griffiths \cite{griffiths1967correlations} and Thompson \cite{thompson1971upper}:
\[
  T_c(J) \leq T^*(J):= \sum_{k}J_{k},\]
where $J_k=J(e_k)$ for the $k$'th standard unit vector $e_k$ in $\Z^d$.

By definition, the Ising model relies on the graph structure of the lattice. A general lattice can be skewed to have different distances between neighbors.
We wish to transport this geometry to the Ising model
by having neighboring spins strongly coupled if they lie close to each other whereas spins far apart should be weakly coupled.
To do so, we will parametrize the space of bonds $J$ using the geometry of a lattice in $\R^d$. 

In principal, we wish to follow a recipe that does this: For any rank $d$ lattice $L$ of $\R^d$, scaled to be of covolume one, take the lengths $\lambda_k$ of the shortest $d$ vectors in $L$ that span $L$. 
Attach an Ising model with bond weights ${J_k}={\lambda_k^{-1}}$ to each lattice. 
Denote the associated critical temperature by $T_c(L)$.
The space of covolume one lattices carries a natural probability measure. 
We shall show: 
\[
	L\mapsto T_c(L)
\]
has finite mean and finite variance.
Viewing the critical temperature as a random variable on the space of lattices enables the application of numerous mathematical tools 
such as the rich representation theory of $\SLR[d]$.
We present evidence supporting this perspective in the literature and provide references in the final section.

\section{Main Result}
\subsection{Ising Model}
The classical nearest neighbor Ising model in $d$-dimensions is defined as follows.
Let $\Lambda$ be a subset of of $\Z^d$,
consider the configurations 
$\sigma\in\Sigma_\Lambda:=\{-1,+1\}^\Lambda$
and the spin random variables
$\sigma_i$ for $i\in \Lambda$
corresponding to the coordinate projections.
The associated Hamiltonian at inverse temperature $\beta$ in a magnetic field of strength $h$ is defined to be
\begin{equation}\label{eq:Hamiltonian}
\mathcal{H}_{J, \beta, \Lambda, h}(\sigma) = -\beta\sum_{i\sim j} J_{ij}\sigma_i\sigma_j - h\sum_{i\in\Lambda}\sigma_i
\end{equation}
where the first sum runs over all tuples of integer vectors $\{i,j\}\in \Lambda\times\Lambda$ that are distance one apart (with respect to the standard graph metric of $\Z^d$).
We assume that the bond variables $J=\{J_{ij}\}$ are strictly positive and translation invariant: $J_{ij}=J(j-i)>0$.
The finite volume Gibbs distribution associated to $\mathcal{H}_{J, \beta,\Lambda, h}$ is the discrete probability measure
\[
  \mu_{J,\beta, \Lambda, h}(\sigma) = \frac1Z_{J, \beta,\Lambda, h}\exp(-\mathcal{H}_{J,\beta,\Lambda, h}(\sigma))
\]
where $Z_{J, \beta,\Lambda, h} = \sum_{\sigma\in\Sigma_{\Lambda}}\exp(-\mathcal{H}_{J,\beta,\Lambda, h}(\sigma))$ is the normalization constant.
Using the magnetization density $m_\Lambda:\Sigma_\Lambda\to[0,1]$, 
$m_\Lambda(\sigma) = \frac{1}{|\Lambda|}\sum_{i\in\Lambda}\sigma_i$,
one defines the spontaneous magnetization as 
\[
  m(J,\beta) = \lim_{h\to 0+}\lim_{\Lambda\to\Z^d}\mu_{J,\beta,\Lambda,h}(m_\Lambda)
\]
where $\Lambda\to\Z^d$ denotes (say) the limit along the sequence $\Lambda_n=\{i\in\Z^d: \|i\|_\infty < n\}$.
The critical inverse temperature $\beta_c(J)$ is then defined to be
\[
\beta_c(J) = \sup\{\beta>0 : m(J,\beta) = 0\}.
\]
It is a fundamental theorem that $\beta_c(J) \in (0,\infty)$ for $d>1$, i.e.\
the Ising ferromagnets exhibit a phase transition at a critical temperature
\[
  T_c(J)=\beta_c(J)^{-1}.
\]
For details, see \cite{friedli2017statistical}.

\subsection{The Geometry of Numbers}
A (full-rank) lattice $L$ in $\R^d$ can be written as $g\Z^d$ for $g\in \GL_d(\R)$, where the column vectors of $g$ define a basis.
The stabilizer group at the standard integer lattice is $\GL_d(\Z)$, and the space of all lattices can be identified with the homogeneous space
$X_{d}:=\GL_d(\R)/\GL_d(\Z)$.
This space captures the geometry of a lattice, in particular how it may go to infinity (leaves eventually any compact set): Letting the covolume go to $0$ or $\infty$,
letting the length of the shortest vector go to zero,
and more generally letting the minimal covolume of a rank $k$ sublattice go to zero.

A lattice $g\Z^d$ is called unimodular if its covolume $\vol(\R^d/g\Z^d)=\det{g}$ is one, and it suffices for us to restrict to such lattices as the degree of freedom for the covolume will already be captured by the temperature for what follows.

$\SL_d(\R)$, the matrix subgroup of $\GL_d(\R)$ consisting of all determinant one matrices, acts transitively on the set of unimodular lattices.
Its stabilizer group is $\SL_d(\Z)$, and the homogeneous space $\SL_d(\R)/\SL_d(\Z)$ defines the space of unimodular lattices $X^1_d$.
It follows from Minkowski's reduction theory that the measure on this quotient induced from the Haar measure of $\SL_d(\R)$ is finite \cite{siegel}.
We normalize it to be probability and denote it $m_{d}$. It is the unique $\SL_d(\R)$-invariant probability measure on $X^1_d$.

For each lattice $L\in X_d^1$ we wish to attach a Hamiltonian $\mathcal{H}_L$, that is, make a choice for $J$ given $L$.
We introduce the Minkowski's successive minima $\lambda_1$,..., $\lambda_d$ of a lattice.
The scalar $\lambda_i$ is defined to be the radius of the smallest closed ball (with respect to some fixed norm $\|\cdot\|$) containing $i$ many linearly independent vectors of $L$,
but any ball of smaller radius contains at most $i-1$ many linearly independent vectors.
It is not true in general, that $\overline{B}_{\lambda_d}$ contains a basis (the first counter-example appears in dimension $d=5$).
However, it is possible to find a basis $v_1$, ..., $v_d$ such that $\lambda_i\leq \|v_i\|\leq 2^{i-1}\lambda_i$, see \cite[Lecture X]{siegel1989lectures}. 
Such basis is called Minkowski reduced, for further details and references see \cite{bremner2011lattice}.

Hence \[J_i=\lambda_i^{-1}\] is an intrinsic choice to the lattice, and will remain to be the definition for the bond weights for the rest of this note.
We selected the inverse norm to reflect that bond weights decrease with distance.
We note that $\lambda_i$ are rotational invariant, that is, for any $k\in\text{O}_d(\R)$, the lattice $L$ and $kL$
share the same Minkowski's successive minima.

As alluded earlier, we note that from this geometric view point, we may think of the temperature $T$ as the
$d$th root of the covolume of $L\in X_d$, as multiplying $L$ by $T$ scales its covolume by $T^d$ and will map $\lambda_i$ to $T\lambda_i$.
Separating therefore the covolume from $L$ (being already captures by the temperature), we have justified the restriction to $L\in X^1_d$.
As noted before, by rotation invariance we placed ourselves on $\text{SO}_d(\R)\backslash X_d^{1}$, the space of shapes of lattices.
For $d=2$, this space is the modular surface $\H/\SLZ$.

In summary, we have obtained a random Ising model by
\[
L=g\Z^d \rightarrow \{\lambda_i=J_i^{-1}\} \rightarrow  \mathcal{H}_{J(L),\beta,\Lambda, 0},
\]
and thus a random variable
\[
L=g\Z^d \mapsto T_c(L).
\]

\subsection{Integrability theorem}

Using reduction theory, the mean field bounds and Onsager's solution for $d=2$, we now state and prove the principal calculation of this note regarding integrability of $L\mapsto T_c(L)$.

\begin{theorem}
  For $L\in X^1_d$ and its Minkowski's successive minima $\{\lambda_i\}$,
  let $T_c(L)$ be the critical temperature associated to the Hamiltonian
$\mathcal{H}_{J,\beta,\Lambda,0}$ where $J_k=\lambda_k^{-1}$.
Let $m_d$ be the Haar probability measure on $X^1_d$.
Then 
\[
	\int_{X^1_d} T_c(L)^p dm_{d}(L)<\infty
\]
for $p<d$ if $d\geq3$ and $p\leq 2$ if $d=2$.
In particular, the average critical temperature $m_{d}(T_c)$ and its variance $\on{var}_{m_{d}}(T_c)$ exist.
\end{theorem}
\begin{proof}

We approximate the integral of $m_{d}$ using Siegel sets, see \cite[Section 1.10]{borel2019introduction}.
  
We start from the Iwasawa decomposition \[\SL_d(\R) =KAN\]
where, $K=\SOR[d]$,
$A$ denotes the $d\times d$-diagonal matrices
and $N$ denote the upper unipotent matrices in $\SL_d(\R)$.

The Haar measure $dg$ of $\SL_d(\R)$ can be decomposed into $dg = \rho(a) dk da dn$, where $dk$, $da$, $dn$ are Haar measures of the groups $K$, $A$, $N$ and $\rho(a) = \prod_{i<j} \frac{a_i}{a_j}$ for $g=kan$, $a=\diag{a_1,\dots,a_d}$.
  
  Define $\mathcal{S}_{t,u}\subset\SL_d(\R)$ by
  $\mathcal{S}_{t,u} = K A_t N_u$, where
\[
    A_t = \{ a = \diag{a_1,\dots,a_d}\; |\; a_i \leq t a_{i+1} \quad (i=1,\dots,d-1) \}
\]
and 
\[
    N_u = \{ n \in N \; |\; |n_{ij}|\leq u \quad (1 \leq i < j \leq d)\}.
\]
Then $\SL_d(\R) = \mathcal{S}_{t,u}\SL_d(\Z)$ if $t\geq \frac{2}{\sqrt{3}}$ and $u\geq \frac12$.
In particular,
\[
  \int_{X^1_d} f(L) dm_{d}(L) \leq \int_{\mathcal{S}_{t,u}} f(g\Z^d) dg\]
for any $f\geq0$.

By the mean field bound,
\[
T_c(g\Z^d)\leq \sum_k J_k(g\Z^d) = \sum_k \lambda_k^{-1}\leq d\lambda_1^{-1} \ll a_1^{-1},
\]
so
\[
m_{X^1_d}(T_c^p)\ll \int_{A_t} a_1^{-p} \rho(a) da.
\]
After a coordinate change, first $b_i=\frac{a_i}{a_{i+1}}$, $\rho(a) = \prod_{i=1}^{d-1} b_i^{i(d-i)}$, then $b_i=\exp y_i$, so
$a_i = b_i\dots b_{d-1}a_d$, $1=a_1\dots a_d = \prod b_i^{i}a_d^{d}$, $a_d = \prod b_i^{-i/d}$, $a_1=\prod b_i^{{(d-i)}/{d}}$
\[
\int_{A_t} a_1^{-p} \rho(a) da = \prod_{1\leq i <d} \int_{-\infty}^{\log t} \exp{\left((i(d-i)-\frac{p(d-i)}{d}) y_i\right)} dy_i
\]
which is finite if $(i-\frac{p}{d})>0$ for $i=1\dots d-1$, i.e. $d>p$

For $d=2$, $T_c$ is known exactly by Onsager's solution \cite{onsager1944crystal} given by the
equation
\[
\sinh(2J_1/T_c)\sinh(2J_2/T_c)=1.
\]
A Taylor approximation (e.g.\ \cite{weng1967critical}) along the limit $J_1/J_2=\eta\to0$ gives an upper bound
\[
T_c(J) \ll \frac{J_1}{-\log\eta}\ll \frac{-1}{\lambda_1\log\lambda_1}. 
\]
This leads to the integral
\[
\int T_c(J)^p \rho(a)da \ll \int_0^t (-a_1\log a_1)^{-p} a_1da_1
\]
which is finite for $p\leq2$. This finishes the theorem also for $d=2$.

\end{proof}

\section{Numerical Integration}
The calculation of $m_2(T_c^p)$ is amenable to numerical tools,
since $T_c$ can be solved by a root-finding procedure from Onsager's solution, and
the quotient
$\SOR[2]\backslash X_2^1=\H/\SL_2(\Z)$ has an easily described
Dirichlet fundamental domain at $i$ given by 
\[
F=\{ z=x+iy\in\H : |z|\geq 1, |x|\leq\frac12\}.\]
The hyperbolic integral is $\int_F f(z) \frac{dxdy}{y^2}$, and we normalize by dividing by the area of $F$, $\pi/3$.
One can use the following recipe to attach a shape of a lattice, $[g]=\SOR[2]g\SLZ[2]$ to a point $z\in F$.
We first note that we may identify $X_2^1$ with the space of lattices up to homothety, $X_2/\R_{>0}=\PSL_2(\R)/\PSL_2(\Z)$.
One goes to $X_2^1$ by picking the unique unimodular representative.
Let $z=x+iy\in F$.
Then $1$,$z$ define a basis of a lattice $L_z$ in $\C\sim\R^2$. We note that $iy$ and $z$ differ by a shear, which is an area preserving transformation.
Hence the covolume of $L_z$ is $y$, and $y^{-\frac12}L_z$ defines a unimodular representative.
The associated basis in $\R^2$ is $(y^{-\frac12},0)^t$ and $(xy^{-\frac12},y^{\frac12})^t$.
We note that  $\lambda_1,\lambda_2$ agree with the length of these vectors, $y^{-\frac12}$, $y^{-\frac12}\sqrt{x^2+y^2}:=y^{-\frac12}r$.

Using the homogeneity property of Onsager's solution, we first solve his implicit formula using a combination of bisection and Newton iterations for $J'=(J_1',J_2')=(1,\frac{J_2'}{J_1'}=\frac{1}{r})$
to obtain $T_c'=T_c(J')$, and then substitute back to get $T_c(\lambda_1^{-1},\lambda_2^{-1}) = T_c'/J_1 = y^{\frac12}T_c(1,1/r)$.

In particular, we only need to solve $T_c$ along a one-dimensional integral after the variable substitution $(x,y)\mapsto (x,r)$, $r^2=x^2+y^2$. In order to integrate  over a finite range, we further substitute $u=1/r$. Numerical integration is done using the extended midpoint method.

The second moment turns out to be slowly convergent, (indeed, it barely converges by the proof of the theorem).
To deal with this instability, we cut off the cusp at some point $r_0$, integrate $F_0 = \{ z\in F: |z|<r_0\}$ as above,
but replace the numerical solution for $T_c(1,1/r)$ by its Taylor approximation
\[
 2\log r \left(1+\frac{\log\log r}{\log r}\right)
\]
inside the cusp $F\setminus F_0$. 
This leads to
\[
m_2(T_c)\approx 2.482 \quad \quad
m_2(T_c^2)\approx 6.979.
\]

\section{Monte Carlo methods}
In higher dimension, we are faced with two difficulties for the estimation of $m_d(T_c)$.
Firstly, the space $X^1_d$ is $d^2-1$ dimensional and explicit fundamental domains are rather cumbersome, see \cite{terras_hecke} surveying the case $d=3$. Secondly, there is no known analytic solution for $T_c$.
Both issues are addressed using Monte Carlo methods.

\subsection{Integration by Sampling}
We replace the integral over $X^1_d$ with the well studied pseudo-random sequence of Hecke points that equidistributes with respect to $m_d$ (\cite{goldstein2003equidistribution}).

Using Hecke points to sample for computer simulations has been proposed in \cite{lubotzky1986hecke}, \cite{lubotzky1987hecke}. 
A recommended introduction is \cite{lubotzky1994discrete}.

Let $p$ be a prime. Then the $p$-Hecke neighbors of a lattice $L$ are defined to be all index $p$ sublattices of $L$. Concretely, if $L=g\Z^d$, then these lattices are given by $gn\Z^d$ where $n$ is one of the following
$1+(p-1)+\dots+(p-1)^{d-1}$ matrices:
\[
  \begin{bmatrix}
    p&   & &  \\
     & 1 & &  \\
     &  & \ddots &  \\
     &   & & 1
    \end{bmatrix}, 
  \begin{bmatrix}
    1& i_1  & &  \\
     & p & &  \\
     &  & \ddots &  \\
     &   & & 1
    \end{bmatrix}, \dots,
  \begin{bmatrix}
    1&   & & i_{1} \\
     & \ddots & & \vdots \\
     &  & 1 & i_{d-1} \\
     &   & & p
    \end{bmatrix} \quad\{i_j =0,\dots p-1\}.
  \]
  One maps them to $X^1_d$ by dividing by $p^{1/d}$.
  Let $\mu_p$ denote the sampling measure over the $p$-Hecke neighbors of $\Z^d$.
  A theoretical bound for the convergence is $|m_d(f)-\mu_p(f)|\ll p^{-\frac{1}{2d^2}}$ given in \cite{goldstein2003equidistribution}.

Applying the method $d=2$ and choosing $p=1336337$ we get 
\( m_2(T_c)\approx 2.486\)  and \( m_2(T_c^2)\approx 6.822\).

\subsection{Temperature Approximation}
For $d=2$, in order to approximate $T_c$ we tried the Swendson-Wang cluster algorithm (\cite{swendsen1987nonuniversal}, \cite{EdSo}) and the invaded cluster algorithm (\cite{machta1995invaded}, \cite{redner2000invaded}) but they did not yield usable results for large ratios of $J_1/J_2$, at least with small spin window approximations on a home computer's CPU.
We got better results using a checkerboard Metropolis algorithm that utilizes the GPU.
We modified a highly optimized CUDA program, which was recently provided by \cite{romero2020high} to allow for anisotropic bond weights. The modification is available here \cite{ising_ani}.
As already observed in \cite{romero2020high}, large system sizes get easily stuck in meta-stable states. 
One of the referees of this note pointed out that \cite{wang2001efficient} and \cite{invernizzi2021opes} might overcome this obstacle.

In the highly anisotropic case this effect amplifies significantly.
By changing the spin window from a square $n\times n$ to a rectangular one, we could mitigate this effect and obtain better results.
Further understanding might allow to choose proper a priori scales for experimentation. 
It might also be worth studying the multiple histogram method for varying $J$ (see e.g.\ \cite{newman1999monte}).

\section{Conclusion and further studies}
The fact that the critical temperature is a function on $L^2(X^1_d)$ has several consequences.
There has been great success in answering classical questions related to the kinetic energy of the Lorentz Gas by moving the view point from one particular lattice to the space of lattices, 
see \cite{marklof2010distribution},\cite{marklof2011boltzmann},\cite{marklof2024kinetic}.
More generally, since $\SLR[d]$ acts on $X^1_d$, one can discuss orbit classifications. 
For example, understanding the orbit (closure) of a particular lattice reveals great insight about the starting point. 
Landmark examples for this approach are \cite{emm}, \cite{lann} and \cite{esmi}.
We provide a list of problems and connections. We hope to inspire the reader to use this note as a starting point to tackle one of these directions.

\subsection{Spectral decomposition of the critical temperature}
The representation theory of $\SLR$ is governed by automorphic forms.
A recent introduction to the topic was published by theoretical physicists \cite{fleig2018eisenstein}. 
What is the spectral decomposition of $T_c$ in $L^2(X^1_2)$? 
Is it continuous, being spanned by Eisenstein series and can essentially be understood by summing over lattice points?
A relevant survey for the lattice point count aspect in Eisenstein series can be found in \cite{burrin2021effective}.
Or does it contain discrete spectrum, which holds additional information about the ambient geometry of the space? 

\subsection{A space of lattices allows for renormalization}
Numerical simulations deep in the cusp have been shown to be difficult. 
Can renormalization ideas, in which a skewed lattice is moved to a generic position be used to overcome these?
We point out that a connection between bond weight anisotropy and the geometry of the finite-size approximate has recently been used for the calculation of the Binder cumulants (\cite{selke2005critical}, \cite{selke2006critical}, \cite{kamieniarz1993universal}). 

\subsection{Dynamics through shears and dilations}
One might study the dynamics of subgroups of the special linear group acting on $X^1_d$.
The action by a diagonal matrix is called dilation, the action by a unipotent matrix a shear.

The full shear has been studied in \cite{ChanLin}, \cite{Hucht}
(a random partial shear is studied in \cite{CiGoSa}) with periodic orbit.
Given a flow speed $v$, with units lattice sites per Monte Carlo sweep
(where an Monte Carlo sweep consists of $n^2$-many Monte Carlo step attempts),
one performs a shear move every $n^2/v$ Monte Carlo steps.
One can imagine to study a shear that equidistributes in $X^1_d$.
Further, one may study a dilation action, for which also closed orbits exist.

\subsection{Limit Approximations}
The limiting approximation of $T_c$ for $d=2$ by means of Taylor approximation has been proven by other methods in \cite{weng1967critical}, \cite{fisher1967critical} that extend to higher dimensions and give the limiting behavior of $T_c$ as
\(
\frac{J_2+\dots+J_d}{J_1}\to0.
\)
The geometry of the space of lattices might guide which other limits to take.
We also left open the question if the integrability exponent $d=p$ is sharp.

\subsection{Quasicrystals}
The anisotropic Ising model on the Penrose tiling, which can be described as a cut-and-project set, has recently be analyzed in \cite{redner2000invaded}.
One can replace the space of lattices by the space of quasicrystals constructed by the cut-and-project method recently introduced
in \cite{msquasi} where a probability measure analogous to $m_d$ is found. 
This allowed them to answer kinetic energy questions on quasicrystals.
Providing bounds for square integrability of a family of counting functions on this space allows to answer questions \textit{in the mean} and from there to deduce them \textit{almost surely}, see \cite{quasischmidt}.

\paragraph{\bf{Acknowledgement}}
This work was supported by the ISF 264/22 while the author was a postdoc at the Weizmann Institute of Science. 
We are grateful for Omri Sarig's hospitality.
We thank Ron Peled and the anonymous referees for helpful comments. 

\bibliographystyle{unsrt}
\bibliography{bib_rene}
\end{document}